\documentclass[11pt]{amsart}
\usepackage{amsmath,amssymb,amsthm}
\usepackage[latin1]{inputenc}
\usepackage{version,tabularx,multicol}
\usepackage{graphicx,float,psfrag}
\usepackage{stmaryrd}

\theoremstyle{plain}
\newtheorem{theo}{Theorem}[section]
 
\newtheorem{cor}[theo]{Corollary}
 
\newtheorem{prop}[theo]{Proposition}
 
\newtheorem{lem}[theo]{Lemma}
 
\newtheorem{conj}[theo]{Conjecture}

\theoremstyle{remark}
\newtheorem{rem}[theo]{Remark}

\newcommand{\mN}{{\mathbb N}}

\newcommand{\mR}{{\mathbb R}}

\newcommand{\mZ}{{\mathbb Z}}

 \def\beqlb{\begin{eqnarray}}\def\eeqlb{\end{eqnarray}}
 \def\beqnn{\begin{eqnarray*}}\def\eeqnn{\end{eqnarray*}}

 \def\ar{\!\!&}

 \def\mbb{\mathbb}

 \def\qed{\hfill$\Box$\medskip}

\newcommand{\bcen}{\begin{center}}
\newcommand{\ecen}{\end{center}}
\newcommand{\bgeqn}{\begin{equation}}
\newcommand{\edeqn}{\end{equation}}

\def\mN{{\mbb  N}}
\def\N{{\mbb  N}}

\def\mR{{\mbb  R}}

\def\mT{{\mbb  T}}
\def\t{{\mathbf t}}
\def\s{{\mathbf s}}

\def\l{\left}
\def\r{\right}

\newcommand{\cal}{\mathcal}
 \def\ar{\!\!\!&}

\begin{document}

\title[Local convergence of critical trees and CB processes]{Local convergence of critical random trees and continuous-state branching processes}
\date{\today}

\author{Xin He}

\address{Xin He, School of Mathematical Sciences, Beijing Normal University, Beijing 100875, P.R.CHINA}

\email{hexin@bnu.edu.cn}

\begin{abstract}
We study the local convergence of critical Galton-Watson trees and L\'{e}vy trees under various conditionings.
Assuming a very general monotonicity property on the functional of random trees,
we show that random trees conditioned to have large functional values always converge locally to immortal trees.
We also derive a very general ratio limit property for functionals of random trees satisfying the monotonicity property.
Then we move on to study the local convergence of critical continuous-state branching processes, and prove a similar result.
Finally we give a definition of continuum condensation trees,
which should be the correct local limits for certain subcritical L\'{e}vy trees under suitable conditionings.
\end{abstract}



\maketitle

\section{Introduction}\label{s:intro}

The local convergence of conditioned Galton-Watson trees (GW trees) has been studied for a long time, dating back to Kesten \cite{K86}, at least.
Over the years, several different conditionings have been studied: large height, large total progeny, and large number of leaves.
Recently, Abraham and Delmas \cite{AD14a,AD14b} provided a convenient framework to study the local convergence of conditioned GW trees, then they used this framework to prove essentially all previous results on the local convergence of conditioned GW trees and also some new ones.
Also very recently, in \cite{H14} we studied the local convergence of GW trees under a new conditioning, which is the conditioning of large maximal out-degree.
An interesting phenomenon is that under any of the conditionings considered in these papers, a conditioned critical GW tree always converges locally to a certain size-biased tree with an infinite spine, which we call an $\emph{immortal tree}$ in this paper. Naturally one would want to ask: Is it true that conditioned critical GW trees always converge locally to immortal trees, under any reasonable conditioning? Is it possible to prove such a general result?

The answer is actually a partial yes. More specifically, we need to distinguish two different formulations of local convergence. We call one formulation the $\emph{tail versions}$ of local convergence, and the other the $\emph{probability versions}$.
For example, let us consider the classical conditioning of large height: If we condition GW trees to have height greater than a large value, then we are considering the tail versions; If we condition GW trees to have height equal to a large value, then the probability versions. For the tail versions, if we assume a very general monotonicity property on the functional of GW trees, then we can prove that critical GW trees conditioned on large functional values always converge locally to immortal trees.
For the probability versions,  it seems less possible to obtain such a general result. Nevertheless, we may impose a more restrictive additivity property on the functional of GW trees and argue with several specific conditionings in mind, to get the probability versions under any of the conditionings that have been studied previously.

Now let us review our results on the local convergence of critical GW trees. In Theorem \ref{GWtail} we prove our general result on the tail versions of local convergence of conditioned critical GW trees. Although this result shows that critical GW trees always converge locally to immortal trees under essentially any conditioning, we only apply it to the conditioning of large width in Corollary \ref{GWwidth}, which is one of our main motivations of this paper. Next we study the corresponding probability versions in Theorem \ref{GWprobability}, where we require a more restrictive additivity property on the functional of GW trees. Then we apply Theorem \ref{GWprobability} to four specific conditionings, which are the conditioning
of large maximal out-degree, the conditioning of large height, the conditioning of large width, and the conditioning of large number
of nodes with out-degree in a given set. Finally we take the argument in the proof of Theorem \ref{GWtail} further to derive in Theorem \ref{DRL} a very general ratio limit property for functionals of GW trees satisfying the monotonicity property. In particular, we give in Proposition \ref{width} two ratio limit results for the width of GW trees.

We have to admit that several results in this paper on the local convergence of conditioned critical GW trees are already known from \cite{AD14a,H14}. Note that a unified method for the local convergence of conditioned critical GW trees has been proposed and used in \cite{AD14a} and also used in \cite{H14} later.
The reason that we revisit all these results here is that we have a different method. Comparing to the method in \cite{AD14a}, we feel that our method has some advantages: First our method seems to be somewhat more direct and intuitive; Second our method seems to be more natural for the proofs of our general results Theorem \ref{GWtail} and Theorem \ref{DRL}; Finally our method can also be used for the local convergence of conditioned critical L\'{e}vy trees, which is also one of our main motivations of this paper. Recall that L\'{e}vy trees are certain scaling limits of GW trees.
Although technically L\'{e}vy trees are more involved than GW trees, we are able to get essentially all the corresponding results for L\'{e}vy trees.

Now let us review our results on the local convergence of conditioned critical L\'{e}vy trees.
Here we only consider the tail versions of local convergence.
Recall that Duquesne \cite{D08} proved the tail versions of local convergence of critical or subcritical L\'{e}vy trees to continuum immortal trees, under the conditioning of large height.
We prove in Theorem \ref{L} that critical L\'{e}vy trees conditioned on large functional values always converge locally to continuum immortal trees,
as long as the functional of L\'{e}vy trees satisfies a very general monotonicity property.
We apply this general result to three specific conditionings, which are the conditioning of large width, the conditioning of large total mass,
and the conditioning of large maximal degree.
Next by taking the argument in the proof of Theorem \ref{L} further,
we derive in Theorem \ref{CRL} a very general ratio limit property for functionals of L\'{e}vy trees satisfying the monotonicity property.
Finally by adapting the proofs of Theorem \ref{L} and Theorem \ref{CRL},
we prove in Theorem \ref{LCCB} that conditioned critical continuous-state branching processes (CB processes) always converge locally to certain CB processes with immigration (CBI processes), again only assuming a very general monotonicity property on the functional of CB processes.

Our method in this paper depends crucially on the criticality of random trees,
so consequently it can not be directly used for the local convergence of subcritical random trees (however, see Corollary \ref{GWfour}).
Recall that relying on the framework and the method of \cite{AD14b}, it has been proved in \cite{H14} that under the conditioning of large maximal out-degree, the local limit of a subcritical GW tree is a condensation tree, which is different from an immortal tree but closely related to it. In the continuous-state setting, it has been shown by Li and He \cite{HL14} that under the conditioning of large maximal jump, the local limit of a subcritical CB process is a certain killed CBI process. Inspired by these two results, we give a precise definition of \emph{continuum condensation trees} in Section \ref{s:CCT}. Naturally we expect continuum condensation trees to be the correct local limits of subcritical L\'{e}vy trees under the conditioning of large maximal degree, however the desired proof seems to be more involved and currently we do not have it yet.
We are confident with this convergence, so we state it explicitly as Conjecture \ref{conjecturemd}.
Under the conditioning of large total progeny,
it is well-known that the local limit of a subcritical GW tree is also a condensation tree, if a certain sum on the offspring distribution is infinite. For a subcritical L\'{e}vy tree, we believe that assuming a similar property on the branching mechanism, the local limit under the conditioning of large total mass is also a continuum condensation tree. We state this convergence explicitly as Conjecture \ref{conjecturetm} and it might be more challenging than Conjecture \ref{conjecturemd}. Finally we consider the conditioning of large width. For both subcritical GW trees and subcritical L\'{e}vy trees, the local convergence under this conditioning is unknown and might be even more challenging than Conjecture \ref{conjecturetm}. So we state it as an open problem to conclude this paper.

This paper is organized as follows.
In Section \ref{s:GW}, we study both the tail versions and the probability versions of local convergence of critical GW trees.
In Section \ref{s:prelim}, we review several basic topics of L\'{e}vy trees.
Section \ref{s:C} is devoted to the local convergence of critical L\'{e}vy trees and CB processes.
Finally in Section \ref{s:CCT}, we define continuum condensation trees and state two conjectures and one open problem related to them.

\section{Local convergence of critical GW trees} \label{s:GW}

In this section first we review several basic topics of GW trees. Then we study the local convergence of critical GW trees, assuming a very general monotonicity property for the tail versions and a more restrictive additivity property for the probability versions. Finally we derive a very general ratio limit property for functionals of GW trees satisfying the monotonicity property.

\subsection{Preliminaries on GW trees} \label{ss:GWprelim}
This section is extracted from \cite{AD14a}. For more details refer to Section 2 in \cite{AD14a}.
Denote by $\mZ_+ = \{0,1,2,\ldots\}$ the set of nonnegative integers and by $\mN = \{1,2,\ldots\}$
the set of positive integers.
Use
$$
\mathcal{U}=\bigcup_{n\geq 0}\mN^n
$$
to denote the set of finite sequences of positive integers with the convention $\mN^0=\{\emptyset\}$.
For $n \geq 1$ and $u = (u_1,\ldots,u_n) \in \mN^n$, let $|u| = n$ be the height of $u$ and $|\emptyset| = 0$ the height of $\emptyset$.
If $u$ and $v$ are two sequences of $\mathcal{U}$, denote by $uv$ the
concatenation of the two sequences,
with the convention that $uv = u$ if $v = \emptyset$  and $uv = v$ if $u = \emptyset$.
The set of ancestors of $u$ is the set
$$
A_u = \{v \in \mathcal{U}: \mbox{ there exists }w \in \mathcal{U}, w \neq \emptyset, \mbox{ such that }u = vw\}.
$$
A tree $\t$ is a subset of $\cal{U}$ that satisfies:
\begin{itemize}
 \item[$\bullet$] $\emptyset \in \t$.

 \item[$\bullet$] If $u\in \t$, then  $A_u\subset \t$.

 \item[$\bullet$] For every $u\in \t$, there exists $k_u(\t) \in \mZ_+$ such that, for every $i\in \mN$, $ui \in \t$ if and only if $1 \leq i \leq k_u(\t)$.
 \end{itemize}
The node $\emptyset$ is called the root of $\t$. The integer $k_u(\t)$ represents the number of offsprings of the node $u$ in the tree $\t$, and we call it the \emph{out-degree} of the node $u$ in the tree $\t$.
The maximal out-degree $M(\t)$ of a tree $\t$ is defined by
\beqlb\label{def:mod}
M(\t)=\sup\{k_u(\t):u\in \t\}.
\eeqlb
The height $H(\t)$ of a tree $\t$ is defined by
\beqlb\label{def:height}
H(\t)=\sup\{|u|:u\in \t\}.
\eeqlb
Denote by $\mT$ the set of trees, by $\mT_0$ the subset of finite trees, and by $\mT^{(h)}$ the subset of trees with height at most $h$,
$$
\mT^{(h)} = \{\t \in \mT: H(\t) \leq h\}.
$$
Finally to get a finite forest of $k$ finite trees, we just put the $k$ finite trees together by keeping all the nodes at their original height, so at height 0 we have $k$ different roots.

For any $\t\in\mT$ and $h \in \mZ_+$,
write $Y_h(\t)$ for the total number of nodes of the tree $\t$ at height $h$.
Also write
$\t_{(h)}=(\t_{(h),i},1\leq i\leq Y_h(\t))$ for the collection of all subtrees above height $h$.
For $h \in \mZ_+$, the restriction function $r_h$ from $\mT$ to $\mT$ is defined by
$$
r_h(\t) =\{ u \in \t: |u| \leq h\}.
$$
We endow the set $\mT$ with the ultra-metric distance
$$
d(\t, \t') = 2^{-\sup \{h\in \mZ_+: \ r_h(\t)=r_h(\t')\}}.
$$
Then a sequence $(\t_n, n\in \mN)$ of trees converges to a tree $\t$ with respect to the distance $d$ if and only if for every $h \in \mN$,
$$
r_h(\t_n) = r_h(\t) \quad \text{for }n\text{ large enough}.
$$
Let $(T_n, n \in \N)$
and $T$ be $\mT$-valued random variables (with respect to the Borel $\sigma$-algebra on $\mT$). Denote by dist $(T)$ the distribution of the
random variable $T$, and denote
$$
\text{dist}(T_n)\rightarrow\text{dist}(T)\quad \text{as }n\rightarrow\infty
$$
for the convergence in distribution of the sequence $(T_n, n \in \N)$ to $T$. It can be proved that the sequence $(T_n, n \in \N)$ converges in distribution to $T$
if and only if for any $h \in \mN$ and $\t\in \mT^{(h)}$,
\beqlb\label{GWconvergence}
\lim_{n\rightarrow \infty}\mathbf{P}[r_h(T_n)=\t]=\mathbf{P}[r_h(T)=\t].
\eeqlb

Let $p=(p_0,p_1,p_2,\ldots)$ be a probability distribution on the set of nonnegative
integers. We exclude the trivial case of $p=(0,1,0,\ldots)$. Denote by $\mu$ the expectation of $p$ and assume that $0<\mu<\infty$.
A $\mT$-valued random variable $\tau$ is a Galton-Watson tree (GW tree) with the offspring distribution
$p$ if the distribution of $k_\emptyset(\tau)$ is $p$ and for $n\in \mN$, conditionally on $\{k_\emptyset(\tau)=n\}$, the subtrees
$(\tau_{(1),1}, \tau_{(1),2},\ldots, \tau_{(1),n})$ are independent and distributed as the original tree $\tau$.
From this definition we can obtain the branching property of GW trees, which says that under the conditional probability $\mathbf{P}[\cdot| Y_h(\tau)= n]$
and conditionally on $r_h(\tau)$,  the subtrees
$(\tau_{(h),1}, \tau_{(h),2},\ldots, \tau_{(h),n})$ are independent and distributed as the original tree $\tau$.
The GW tree is called critical (resp. subcritical, supercritical) if $\mu = 1$ (resp. $\mu < 1$,
$\mu >1$). In the critical and subcritical case, we have that a.s. $\tau$ belongs to $\mT_0$.

Immortal trees can be defined for critical or subcritical offspring distributions.
We recall the following definition from Section 1 in \cite{AD14b}, which first appeared in Section 5 of \cite{J12}. Let $p$ be a critical or subcritical offspring distribution.
Let $\tau^*(p)$ denote the random
tree which is defined by:
\begin{itemize}
 \item[i)] There are two types of nodes: $normal$ and $special$.

 \item[ii)] The root is special.

 \item[iii)] Normal nodes have offspring distribution $p$.

 \item[iv)] Special nodes have offspring distribution the size-biased distribution $\hat{p}$ on $\mZ_+$
defined by $\hat{p}_k =kp_k/\mu$ for $k \in \mZ_+$.

  \item[v)] The offsprings of all the nodes are independent of each others.

   \item[vi)] All the children of a normal node are normal.

    \item[vii)] When a special node gets several children, one of them is selected uniformly
at random and is special while the others are normal.
 \end{itemize}
Notice that a.s. $\tau^*(p)$ has one unique infinite spine. We call it an $\emph{immortal tree}$. By the definitions of GW trees and immortal trees, it can be shown that
for any $b \in \mZ_+$ and $\t\in \mT^{(b)}$,
\beqlb\label{Kesten's tree}
\mu^b\,\mathbf{E}[\mathbf{1}\{r_b(\tau^*(p))=\t\}] =\mathbf{E}[Y_b(\tau)\mathbf{1}\{r_b(\tau)=\t\}].
\eeqlb

\subsection{The tail versions of local convergence}\label{ss:GWtail}

Let $A$ be a nonnegative integer-valued function defined on $\mT_0$.
Recall that for any $\t\in\mT_0$, we write $(\t_{(b),i},1\leq i\leq Y_b(\t))$ for the collection of all subtrees above height $b$.
We introduce the following monotonicity property of $A$:
\beqlb \label{GWmonotonicity}
A(\t_{(b),i})\leq A(\t), \quad \text{for any } b\in \mN \text{ and } 1\leq i\leq Y_b(\t).
\eeqlb
The meaning of this monotonicity property (\ref{GWmonotonicity}) should be clear: For any $b\in \mN$, the value of $A$ on the whole tree is not less than that on any subtree above height $b$.

Define $v_n=\mathbf{P}[A (\tau)>n]\in [0,1]$ and $\mathbf{P}_n[\cdot]=\mathbf{P}[\cdot|A (\tau)>n]$ when $v_n>0$.
The following theorem asserts that if the monotonicity property (\ref{GWmonotonicity}) holds for $A$, then
under the conditional probability $\mathbf{P}_n$, the GW tree $\tau(p)$ converges locally to the immortal tree $\tau^*(p)$.

\begin{theo} \label{GWtail}
Assume that $p$ is critical and $v_n>0$ for all $n$.
If $A$ satisfies the monotonicity property (\ref{GWmonotonicity}),
then as $n\rightarrow\infty$,
$$
\normalfont\text{dist }(\tau\big|A(\tau)>n)\rightarrow\text{dist }(\tau^*).
$$
\end{theo}

\proof By (\ref{GWconvergence}), we only have to prove that for any $b \in \mN$ and $\t\in \mT^{(b)}$,
$$
\lim_{n\rightarrow\infty} \mathbf{P}_n[r_b(\tau)=\t]
=\mathbf{P}[r_b(\tau^*)=\t].
$$
Recall from (\ref{Kesten's tree}) that when $\mu=1$, for any $b \in \mN$ and $\t\in \mT^{(b)}$,
$$
\mathbf{P}[r_b(\tau^*)=\t]
=\mathbf{E}[\mathbf{1}\{r_b(\tau^*)=\t\}]
=\mathbf{E}[Y_b(\tau) \mathbf{1}\{r_b(\tau)=\t\}].
$$
So it suffices to show that for any $b \in \mN$ and $\t\in \mT^{(b)}$,
\beqlb \label{GWenough}
\lim_{n\rightarrow\infty} \mathbf{E}_n[\mathbf{1}\{r_b(\tau)=\t\}]
=\mathbf{E}[Y_b(\tau) \mathbf{1}\{r_b(\tau)=\t\}].
\eeqlb

To prove (\ref{GWenough}), first recall that if the value of $A$ on a subtree above height $b$ is greater than $n$, than the value of $A$ on the whole tree is greater than $n$, by the monotonicity property (\ref{GWmonotonicity}). Then recall from the branching property in Section \ref{ss:GWprelim} that
under $\mathbf{P}$ and conditional on $r_b(\tau)$ the probability that the value of $A$ is greater than $n$ on at least one subtree above height $b$ is
$$1-(1- v_n)^{Y_b(\tau)}.$$
So the monotonicity property and the branching property imply that
$$
\mathbf{E}\l[\mathbf{1}\{A(\tau)>n\}\mathbf{1}\{r_b(\tau)=\t\}\r]\geq
\mathbf{E}\l[\l(1-(1- v_n)^{Y_b(\tau)}\r)\mathbf{1}\{r_b(\tau)=\t\}\r].
$$
Thus we see that as $n\rightarrow\infty$,
\beqnn
\mathbf{E}_n\l[\mathbf{1}\{r_b(\tau)=\t\}\r]
\ar\geq\ar\frac{1}{v_n}\mathbf{E}\l[\l(1-(1- v_n)^{Y_b(\tau)}\r)\mathbf{1}\{r_b(\tau)=\t\}\r] \cr
\ar\rightarrow\ar \mathbf{E}\l[Y_b(\tau)\mathbf{1}\{r_b(\tau)=\t\}\r],\cr
\eeqnn
where the convergence follows from the monotone convergence. Note that $v_n\rightarrow 0$ as $n\rightarrow \infty$.

From the previous paragraph we get the inequality that
$$
\liminf_{n\rightarrow\infty} \mathbf{E}_n\l[\mathbf{1}\{r_b(\tau)=\t\}\r]
\geq \mathbf{E}\l[Y_b(\tau)\mathbf{1}\{r_b(\tau)=\t\}\r].
$$
Note that
$$
\sum_{\t\in\mT^{(b)}}\mathbf{E}_n\l[\mathbf{1}\{r_b(\tau)=\t\}\r]=1
\quad\text{and}\quad
\sum_{\t\in\mT^{(b)}}\mathbf{E}\l[Y_b(\tau)\mathbf{1}\{r_b(\tau)=\t\}\r]=1,
$$
since $\mathbf{E}[Y_b(\tau)]=1$ by (\ref{Kesten's tree}).
This implies that all inequalities above are actually equalities, so we have proved (\ref{GWenough}).
\qed

Note that it is easy to think of a conditioning under which the local limits of conditioned critical GW trees are not immortal trees, such as the conditioning of large minimal out-degree, where the minimal out-degree of a tree is defined to be the minimum of positive out-degrees of all nodes in the tree.
It should be clear that the minimal out-degree does not satisfy the monotonicity property (\ref{GWmonotonicity}).

Although Theorem \ref{GWtail} holds for any $A$ satisfying the monotonicity property (\ref{GWmonotonicity}), one of our original motivations for this result is the local convergence under the conditioning of large width. So right now we will only apply Theorem \ref{GWtail} to this specific conditioning. Here the width $W(\t)$ of a tree $\t$ is defined to be $\sup_{b\in\mZ_+} Y_b(\t)$. Note that $\mathbf{P}[W(\tau)>n]>0$ for any $n$ if and only if $p_0+p_1<1$. Recall that we exclude the trivial case of $p_1=1$. Then Theorem \ref{GWtail} immediately gives the local convergence of critical GW trees to immortal trees, under the conditioning of large width.

\begin{cor} \label{GWwidth}
Assume that $\mu=1$.
Then as $n\rightarrow\infty$,
$$
\normalfont\text{dist }(\tau\big|W(\tau)>n)\rightarrow\text{dist }(\tau^*).
$$
\end{cor}

\subsection{The probability versions of local convergence}\label{ss:GWprobability}
The probability versions automatically imply the corresponding tail versions,
since the tail versions can be written as sums of the corresponding probability versions.
More precisely, we have
$$
\mathbf{P}[\cdot|A (\tau)>n]=\sum_{m>n}\mathbf{P}[\cdot|A (\tau)=m]\frac{\mathbf{P}[A (\tau)=m]}{\mathbf{P}[A (\tau)>n]}.
$$
To get the probability versions,
we have to impose a more restrictive additivity property on the functional $A$, which is similar in spirit to the additivity property (3.1) in \cite{AD14a}.

Let $A$ be a nonnegative integer-valued function defined on the space of finite forests.
Recall that $\t_{(b)}$ is the sub-forest of the tree $\t$ above height $b$ and $r_b(\t)$ is the subtree of the tree $\t$ below height $b$.
We introduce the following additivity property of $A$: For any fixed $b\in \mN$ and $\s\in \mT^{(b)}$,
\beqlb \label{GWadditivity}
A(\t)=A(\t_{(b)})+B(r_b(\t)), \quad \text{for large enough }A(\t) \text{ with }r_b(\t)=\s,
\eeqlb
where $B$ is a nonnegative integer-valued function on $\mT_0$.

Define $v_{(n)}=\mathbf{P}[A (\tau)=n]\in [0,1]$ and $\mathbf{P}_{(n)}[\cdot]=\mathbf{P}[\cdot|A (\tau)=n]$ when $v_{(n)}>0$.
Let $\tau^{(k)} = (\tau_1,\ldots,\tau_k)$ be the forest of $k$ i.i.d. GW trees with offspring distribution $p$. Write $v_{(n)}(k)=\mathbf{P}[A (\tau^{(k)})=n]$.
The following Theorem asserts that if the additivity property (\ref{GWadditivity}) holds for $A$ and some additional ratio limit properties hold for $v_{(n)}$ and $v_{(n)}(k)$, then
under the conditional probability $\mathbf{P}_{(n)}$, the GW tree $\tau(p)$ also converges locally to the immortal tree $\tau^*(p)$.

\begin{theo} \label{GWprobability}
Assume that $v_n>0$ for all $n$, the additivity property (\ref{GWadditivity}) holds for $A$,
\beqlb \label{GWk}
\liminf_{n\rightarrow\infty} v_{(n)}(k)/v_{(n)}\geq k,\quad \text{for any }k\in\mN,
\eeqlb
and one of the following two conditions holds:
\begin{itemize}
 \item[I.] $\mu=1$, and $\limsup_{n\rightarrow\infty}v_{(n)}/v_{(n-B(\t))}\leq 1$ for any $\t\in \mT_0$.

\item[II.] $0<\mu\leq 1$, $B(\t)=H(\t)$  for $\t\in \mT_0$, and $\limsup_{n\rightarrow\infty}v_{(n+1)}/v_{(n)}\leq \mu$.
\end{itemize}
Then as $n\rightarrow\infty$,
$$
\normalfont\text{dist }(\tau\big|A(\tau)=n)\rightarrow\text{dist }(\tau^*).
$$
Note that all limits are understood along the infinite sub-sequence $\{n: v_{(n)} > 0\}$.
\end{theo}

\proof For Case I, by the additivity property (\ref{GWadditivity}) and the branching property, we see that
for any $b\in \mN$ and $\t\in \mT^{(b)}$, when $n$ is large enough,
\beqnn
\mathbf{E}_{(n)}\l[\mathbf{1}\{r_b(\tau)=\t\}\r]
\ar=\ar\frac{1}{v_{(n)}}\mathbf{E}\l[\mathbf{1}\{A(\tau_b)=n-B(r_b(\tau))\}\mathbf{1}\{r_b(\tau)=\t\}\r]\cr
\ar=\ar\frac{1}{v_{(n)}}\mathbf{E}\l[v_{(n-B(r_b(\tau)))}(Y_b(\tau))\mathbf{1}\{r_b(\tau)=\t\}\r].\cr
\eeqnn
Then by our assumptions and the Fatou's lemma, we get
$$
\liminf_{n\rightarrow\infty}\mathbf{E}_{(n)}\l[\mathbf{1}\{r_b(\tau)=\t\}\r]
\geq\mathbf{E}\l[Y_b(\tau)\mathbf{1}\{r_b(\tau)=\t\}\r].
$$
From the first paragraph and the last paragraph of the proof of Theorem \ref{GWtail}, we see that the above inequality is enough to imply the local convergence for Case I.

The proof of Case II is similar. We first argue that
\beqnn
\liminf_{n\rightarrow\infty}\mathbf{E}_{(n)}\l[\mathbf{1}\{r_b(\tau)=\t\}\r]
\ar=\ar\liminf_{n\rightarrow\infty}\frac{1}{v_{(n)}}\mathbf{E}\l[\mathbf{1}\{A(\tau_b)=n-b\}\mathbf{1}\{r_b(\tau)=\t\}\r] \cr
\ar=\ar\liminf_{n\rightarrow\infty}\frac{1}{v_{(n)}}\mathbf{E}\l[v_{(n-b)}(Y_b(\tau))\mathbf{1}\{r_b(\tau)=\t\}\r] \cr
\ar\geq\ar \mu^{-b}\mathbf{E}\l[Y_b(\tau)\mathbf{1}\{r_b(\tau)=\t\}\r].\cr
\eeqnn
Since $\mu^{-b}\mathbf{E}[Y_b(\tau)]=1$, clearly the above inequality is also enough to imply the local convergence for Case II.
\qed

Next we will apply Theorem \ref{GWprobability} to four specific conditionings, which are
the conditioning of large height,
the conditioning of large maximal out-degree,
the conditioning of large width,
and the conditioning of large number of nodes with out-degree in a given set.
First we show in the following lemma that Condition (\ref{GWk}) in
Theorem \ref{GWprobability} holds when a certain maximum property holds for the functional $A$.
This result will be applied to the maximal out-degree and the height.

\begin{lem} \label{GWmax}
Assume that $A(\t_{(b)})=\max_{1\leq i\leq  Y_b(\t)} A(\t_{(b),i})$ for any $\t\in \mT_0$ and $b\in \mN$, and $v_n>0$ for all $n$.
Then for any $k\in \mN$,
$$
\lim_{n\rightarrow\infty} v_{(n)}(k)/v_{(n)}= k,
$$
where the limit is understood along the infinite sub-sequence $\{n: v_{(n)} > 0\}$.
\end{lem}

\proof Just notice that
$$
v_{(n)}(k)=(1- v_n)^k-(1- v_{n-1})^k=v_{(n)}\l[\sum_{0\leq i \leq k-1} (1- v_n)^{k-1-i}(1- v_{n-1})^i\r].
$$
Then since $v_n\rightarrow 0$ as $n\rightarrow\infty$,
$$
\lim_{n\rightarrow\infty}\sum_{0\leq i \leq k-1} (1- v_n)^{k-1-i}(1- v_{n-1})^i= k.
$$
\qed

For the forest $\t^{(k)}=(\t_1,\t_2,\ldots,\t_k)$, again we write $Y_h(\t^{(k)})$ for the total number of nodes in the forest $\t^{(k)}$ at height $h$, that is, $Y_h(\t^{(k)})=\sum_{1\leq i\leq k}Y_h(\t_i)$.
Then define $W(\t^{(k)})=\sup_{h\in\mZ_+}Y_h(\t^{(k)})$.

\begin{lem} \label{GWwidth>}
Assume that $A=W$, and the critical or subcritical $p$ has bounded support with $p_0+p_1<1$.
Then for any $k\in \mN$,
$$
\liminf_{n\rightarrow\infty} v_{(n)}(k)/v_{(n)}\geq k,
$$
where the limit is understood along the infinite sub-sequence $\{n: v_{(n)} > 0\}$.
\end{lem}

\proof Assume that $k\geq 2$. Let $N=\sup\{n: p_n>0\}<\infty$ be the supremum of the support of $p$. Use $\lfloor r \rfloor$ to denote the largest integer less than or equal to $r$.
We argue that if $W(\tau_1)=n$ and $H(\tau_i)<\lfloor n/(kN)\rfloor$ for $2\leq i \leq k$, then the width of $\tau^{(k)}$ strictly below generation $\lfloor n/(kN)\rfloor$ is less than
$\lfloor n/(kN)\rfloor kN \leq n$, that is,
$$
\sup_{h< \lfloor n/(kN)\rfloor}Y_h(\tau^{(k)})< \lfloor n/(kN)\rfloor kN \leq n,
$$
which implies that $W(\tau_1)=n$ is achieved after generation $\lfloor n/(kN)\rfloor$ and $W(\tau^{(k)})=W(\tau_1)=n$.
Using this observation, we see that
\beqnn
\liminf_{n\rightarrow\infty} \frac{v_{(n)}(k)}{v_{(n)}}\ar\geq \ar\liminf_{n\rightarrow\infty}\frac{k\mathbf{P}[W(\tau_1)=n,H(\tau_i)<\lfloor n/(kN)\rfloor,2\leq i \leq k]}{\mathbf{P}[W(\tau)=n]}\cr
\ar=\ar\liminf_{n\rightarrow\infty}k\l(\mathbf{P}\l[H(\tau)<\lfloor n/(kN)\rfloor\r]\r)^{k-1}\cr
\ar=\ar k.
\eeqnn
\qed

We turn to the conditioning of large number of nodes with out-degree in a given set, which is the main topic of \cite{AD14a,AD14b}.
For any $\mathcal{A}\subset\mZ_+$, denote by $L_{\mathcal{A}}(\t)$ the total number of nodes in the tree $\t$ with out-degree in $\mathcal{A}$. For example, $L_{\mZ_+}(\t)$ is just the total progeny of the tree $\t$, and $L_{\{0\}}(\t)$ is just the total number of leaves of the tree $\t$.

Now we show that when combined with several results from \cite{AD14a} (which are not directly related to the local convergence), our Theorem \ref{GWprobability} can also be used to prove all the known probability versions of local convergence of critical GW trees from \cite{AD14a,H14}. Recall (\ref{def:mod}) and (\ref{def:height}), the definitions of the maximal out-degrees and the height of trees. For $\mathcal{A}\subset\mZ_+$, define $p(\mathcal{A})=\sum_{k\in\mathcal{A}}p_k$.

\begin{cor} \label{GWfour}
Take any $\mathcal{A}\subset\mZ_+$ with $p(\mathcal{A})>0$.
If $p$ is critical, then as $n\rightarrow\infty$,
$$
\normalfont
\normalfont\text{dist }(\tau\big|L_{\mathcal{A}}(\tau)=n)\rightarrow\text{dist }(\tau^*)
\quad\text{and}\quad
\text{dist }(\tau\big|M(\tau)=n)\rightarrow\text{dist }(\tau^*),
$$
where the limits are understood along the infinite sub-sequences $\{n\in\mN: \mathbf{P}(L_{\mathcal{A}}(\tau)=n)>0\}$ and $\{n\in\mN: p_n>0\}$, respectively.
If $p$ is critical with bounded support, then as $n\rightarrow\infty$,
$$
\normalfont\text{dist }(\tau\big|W(\tau)=n)\rightarrow\text{dist }(\tau^*),
$$
where the limit is understood along the infinite sub-sequence $\{n\in\mN: \mathbf{P}(W(\tau)=n)>0\}$.
If $p$ is critical or subcritical, then as $n\rightarrow\infty$,
$$
\normalfont\text{dist }(\tau\big|H(\tau)=n)\rightarrow\text{dist }(\tau^*).
$$
\end{cor}

\proof For the conditioning of large maximal out-degree, clearly we may let $A(\t)=M(\t)$, then let $A(\t_{(b)})=\max_{1\leq i\leq Y_b(\t)} A(\t_{(b),i})$ and $B\equiv0$.
Now the local convergence follows from Case I in Theorem \ref{GWprobability}, Lemma \ref{GWmax}, and the simple fact that for any $n>0$, $\mathbf{P}[M(\tau)=n]>0$ if and only if $p_n>0$.

For the conditioning of large width, we let $A(\t)=W(\t)=\sup_a Y_a(\t)$, $A(\t_{(b)})=\sup_a Y_a(\t_{(b)})$, and $B\equiv0$.
Now the local convergence follows from Case I in Theorem \ref{GWprobability}, Lemma \ref{GWwidth>}, and the simple fact that $\mathbf{P}[W(\tau)>n]>0$ for any $n$ if and only if $p_0+p_1<1$.

For the conditioning of large height, we let $A(\t)=H(\t)$, $A(\t_{(b)})=\max_{1\leq i\leq Y_b(\t)} A(\t_{(b),i})$ and $B(\t)=H(\t)$.
Now the local convergence follows from (4.5) in \cite{AD14a}, Case II in Theorem \ref{GWprobability}, Lemma \ref{GWmax}, and the trivial fact that $\mathbf{P}[H(\tau)=n]>0$ for any $n\in \mN$.

For the conditioning of large number of nodes with out-degree in a given set $\mathcal{A}\subset\mZ_+$, we only give an outline of our proof and leave the details to the reader.
let $A(\t)=L_{\mathcal{A}}(\t)$, then let $A(\t_{(b)})=\sum_{1\leq i\leq Y_b(\t)} A(\t_{(b),i})$ and $B(\t)=L_{\mathcal{A}}(\t)-Y_{H(\t)}(\t)\mathbf{1}\{0\in \mathcal{A}\}$.
First consider the case of $\mathcal{A}=\mZ_+$.
Then the local convergence follows from Theorem \ref{GWprobability}, the Dwass formula, and a strong ratio theorem for random walks. See e.g., (4.3) in \cite{AD14b} for the Dwass formula and (8.2) in \cite{AD14b} for the strong ratio theorem. These two results combined imply Condition (\ref{GWk}) and Condition I in our Theorem \ref{GWprobability}.
Finally by Section 5.1 in \cite{AD14a}, we know that the case of any general $\mathcal{A}\subset\mZ_+$ can be reduced to the case of $\mathcal{A}=\mZ_+$
in the sense that $L_{\mathcal{A}}$ of any critical GW tree equals $L_{\mZ_+}$ of a corresponding critical GW tree.
So for any $\mathcal{A}\subset\mZ_+$, Condition (\ref{GWk}) and Condition I in our Theorem \ref{GWprobability} hold since they hold for $\mZ_+$.

\qed

\subsection{A general ratio limit property}\label{ss:DRL}

Let $A$ be a nonnegative integer-valued function defined on the space of finite forests.
We introduce the following monotonicity property of $A$:
\beqlb \label{Dmonotonicity}
A(\t_{(b),i})\leq A(\t_{(b)})\leq A(\t), \quad \text{for any } \t\in \mT_0,\ b\in \mN,  \text{ and } 1\leq i \leq Y_b(\t).
\eeqlb

Write $\mathbf{P}[A>n]$ for $\mathbf{P}[A (\tau)>n]$ and $\mathbf{P}^{(k)}[A>n]$ for $\mathbf{P}[A (\tau^{(k)})>n]$.
The following theorem asserts that if the monotonicity property (\ref{Dmonotonicity}) holds for $A$, then as $n\rightarrow\infty$,
$\mathbf{P}^{(k)}[A>n]$ is always asymptotically equivalent to $k\mathbf{P}[A>n]$.

\begin{theo} \label{DRL}
Assume that $p$ is critical, $\mathbf{P}[A>n]>0$ for all $n$, and $A$ satisfies the monotonicity property (\ref{Dmonotonicity}).
Then for any $k\in \mN$,
$$
\lim_{n\rightarrow\infty} \frac{\mathbf{P}^{(k)}[A>n]}{\mathbf{P}[A>n]}=k.
$$
Assume additionally that for some $b\in\mN$, $\s\in\mT^{(b)}$ with $\mathbf{P}[r_b(\tau)=\s]>0$, and $r>0$,
$A(\t)-A(\t_{(b)})=r$ for large enough $A(\t)$ with $r_b(\t)=\s$.
Then for any $k\in \mN$ and $r\in \mN$,
$$
\lim_{n\rightarrow\infty} \frac{\mathbf{P}^{(k)}[A>n-r]}{\mathbf{P}^{(k)}[A>n]}=1.
$$
\end{theo}

\proof First as in the proof of Theorem \ref{GWtail}, for any $k\in\mN$,
\beqlb \label{D>}
\liminf_{n\rightarrow\infty} \frac{\mathbf{P}^{(k)}[A>n]}{\mathbf{P}[A>n]}\geq \liminf_{n\rightarrow\infty} \frac{1-\l(1-\mathbf{P}[A>n]\r)^k}{\mathbf{P}[A>n]}= k.
\eeqlb

Next we argue that for any $k\in\mN$, if there exists some $b\in\mN$ with $\mathbf{P}[Y_b(\tau)=k]>0$, then
\beqlb \label{D=}
\lim_{n\rightarrow\infty} \frac{\mathbf{P}^{(k)}[A>n]}{\mathbf{P}[A>n]}= k.
\eeqlb
To prove this, pick $\t$ with $\mathbf{P}[r_b(\tau)=\t]>0$ and $Y_b(\t)=k$.
As in the proof of Theorem \ref{GWtail}, we have
$$
\lim_{n\rightarrow\infty}\mathbf{P}[r_b(\tau)=\t|A>n]=\lim_{n\rightarrow\infty}\frac{\mathbf{P}^{(k)}[A>n]}{\mathbf{P}[A>n]}\mathbf{P}[r_b(\tau)=\t]=k\mathbf{P}[r_b(\tau)=\t],
$$
which implies (\ref{D=}).

Finally assume that for some $k\in\mN$,
$$
\limsup_{n\rightarrow\infty} \frac{\mathbf{P}^{(k)}[A>n]}{\mathbf{P}[A>n]}>k.
$$
By the facts that $\mathbf{E}[Y_b(\tau)]=1$ and a.s. $\lim_{b\rightarrow\infty}Y_b(\tau)=0$, we can pick some $k'\in\mN$ such that $k'>k$
and there exists some $b\in\mN$ with $\mathbf{P}[Y_b(\tau)=k']>0$.
So (\ref{D=}) holds for $k'$. However as in (\ref{D>}), we also have
\beqnn
\mathbf{P}^{(k')}[A>n]\ar\geq\ar 1-(1-\mathbf{P}^{(k)}[A>n])(1-\mathbf{P}[A>n])^{k'-k}\cr
\ar=\ar 1-(1-\mathbf{P}[A>n])^{k'-k}+\mathbf{P}^{(k)}[A>n](1-\mathbf{P}[A>n])^{k'-k},\cr
\eeqnn
which implies that
$$
\limsup_{n\rightarrow\infty} \frac{\mathbf{P}^{(k')}[A>n]}{\mathbf{P}[A>n]}
\geq (k'-k)+\limsup_{n\rightarrow\infty}\frac{\mathbf{P}^{(k)}[A>n]}{\mathbf{P}[A>n]}>k',
$$
a contradiction to (\ref{D=}) for $k'$.

For the second statement, by the assumptions and the argument in the proof of Theorem \ref{GWtail}, we have for $k=Y_b(\s)$ and the particular $r$ in the assumptions,
$$
\lim_{n\rightarrow\infty}\mathbf{P}[r_b(\tau)=\s|A>n]=\lim_{n\rightarrow\infty}\frac{\mathbf{P}^{(k)}[A>n-r]}{\mathbf{P}[A>n]}\mathbf{P}[r_b(\tau)=\s]=k\mathbf{P}[r_b(\tau)=\s],
$$
which implies that for $k=Y_b(\s)$ and this particular $r$,
$$
\lim_{n\rightarrow\infty} \frac{\mathbf{P}^{(k)}[A>n-r]}{\mathbf{P}[A>n]}=k.
$$
Then by the first statement we have for this particular $r$,
$$
\lim_{n\rightarrow\infty} \frac{\mathbf{P}[A>n-r]}{\mathbf{P}[A>n]}=1,
$$
which implies that for any $r\in\mN$,
$$
\lim_{n\rightarrow\infty} \frac{\mathbf{P}[A>n-r]}{\mathbf{P}[A>n]}=1.
$$
Finally apply the first statement again to finish the proof of the second statement.
\qed

For the probability $\mathbf{P}[A=n]$, it seems not possible to obtain a general result like Theorem \ref{DRL}. However for some specific functional $A$, it is possible to improve the inequality (\ref{GWk}) in Theorem \ref{GWprobability} to an equality.
Recall Theorem 1 in \cite{L76} and Theorem 1 in \cite{BV96}. The following proposition might be regarded as a generalization of these two results.

\begin{prop} \label{width}
Assume that $p$ is critical.
Then for any $k\in \mN$,
$$
\lim_{n\rightarrow\infty} \frac{\mathbf{P}^{(k)}[W>n]}{\mathbf{P}[W>n]}=k.
$$
Assume additionally that $p$ has bounded support.
Then for any $k\in \mN$,
$$
\lim_{n\rightarrow\infty} \frac{\mathbf{P}^{(k)}[W=n]}{\mathbf{P}[W=n]}=k,
$$
where the limit is understood along the infinite sub-sequence $\{n: \mathbf{P}[W=n] > 0\}$.
\end{prop}

\proof The first statement is immediate from Theorem \ref{DRL}.

For the second statement, first recall Lemma \ref{GWwidth>} and the fact that $W$ satisfies the additivity property (\ref{GWadditivity}) with $B\equiv0$.
As in the proof of Theorem \ref{DRL}, we have for any $k\in\mN$, if there exists some $b\in\mN$ with $\mathbf{P}[Y_b(\tau)=k]>0$, then
\beqlb \label{D==}
\lim_{n\rightarrow\infty} \frac{\mathbf{P}^{(k)}[W=n]}{\mathbf{P}[W=n]}= k.
\eeqlb
Now assume that for some $k\in\mN$,
$$
\limsup_{n\rightarrow\infty} \frac{\mathbf{P}^{(k)}[W=n]}{\mathbf{P}[W=n]}>k.
$$
As in the proof of Theorem \ref{DRL}, there exists some $k'\in\mN$ such that $k'>k$ and (\ref{D==}) holds for $k'$. However as in the proof of Lemma \ref{GWwidth>}, we also have
\beqnn
\mathbf{P}^{(k')}[W=n]\geq \ar\ar\mathbf{P}^{(k)}[W=n]\l(\mathbf{P}\l[H<\lfloor n/(k'N)\rfloor\r]\r)^{k'-k}\cr
\ar\ar+(k'-k)\mathbf{P}[W=n]\l(\mathbf{P}\l[H<\lfloor n/(k'N)\rfloor\r]\r)^{k'-1},\cr
\eeqnn
which implies that
$$
\limsup_{n\rightarrow\infty} \frac{\mathbf{P}^{(k')}[W=n]}{\mathbf{P}[W=n]}
> k+(k'-k)=k',
$$
a contradiction to (\ref{D==}) for $k'$.
\qed

\section{Preliminaries on L\'{e}vy trees}\label{s:prelim}

This section is extracted from \cite{D08}. For more details refer to Section 1.2, 3.1, and 3.3 in \cite{D08}.

\subsection{Branching mechanisms of L\'{e}vy trees}\label{ss:BM}

We consider a L\'{e}vy tree with the branching mechanism
\beqlb \label{BM}
\Phi(\lambda)=\alpha \lambda+\beta \lambda^2+\int_{(0,\infty)}\pi(d\theta)(e^{-\lambda \theta}-1+\lambda \theta),
\eeqlb
where $\alpha\geq 0$, $\beta\geq 0$, and the \emph{L\'{e}vy measure} $\pi$ is a $\sigma$-finite measure on $(0,\infty)$ satisfying
$\int_{(0,\infty)}\pi(d\theta)(\theta \wedge \theta^2)<\infty$.
When we talk about height processes of L\'{e}vy trees (see Section \ref{ss:bp}), we always assume the condition
\beqlb \label{assumptionstrong}
\int_1^\infty 1/\Phi(\lambda)d\lambda <\infty,
\eeqlb
which implies that
\beqlb \label{assumptionweak}
\beta>0 \quad \text{or} \quad \int_{(0,1)}\theta\pi(d\theta)=\infty.
\eeqlb

We then consider a spectrally positive L\'{e}vy process $X=(X_t,t\geq 0)$ with the Laplace exponent $-\Phi$, that is, for $\lambda\geq 0$,
\beqlb \label{Levy}
\mathbf{E}[\exp(-\lambda X_t)]=\exp[t \Phi(\lambda)].
\eeqlb
We also consider a bivariate subordinator $(U,V)=((U_t,V_t),t\geq 0)$, that is, a
$[0,\infty)\times[0,\infty)$-valued L\'{e}vy process started at $(0,0)$ (see e.g., Page 162 in \cite{K14}). Its distribution is characterized by the
Laplace exponent $\Phi(p,q)$: For $p,q\geq 0$,
\beqlb \label{subordinator}
\mathbf{E}[\exp(-p U_t-q V_t)]=\exp[-t \Phi(p,q)],
\eeqlb
where
$$
\Phi(p,q)=\frac{\Phi(p)-\Phi(q)}{p-q}-\alpha \quad\text{for} \quad p\neq q, \quad\text{and}\quad \Phi(p,p)=\Phi'(p)-\alpha.
$$

\subsection{The excursion representation of CB processes}\label{ss:CB}

We also consider a continuous-state branching process (CB process) $Y=(Y_t,t\geq 0)$ with the branching mechanism $\Phi$ given in (\ref{BM}).
The branching mechanism and the corresponding CB process and L\'{e}vy tree are called \emph{subcritical} if $\alpha>0$ and \emph{critical} if $\alpha=0$.
Let $\mathbf{P}_x[Y\in\cdot]$ be the distribution of $Y$ under the assumption of $Y_0=x$, and $\mathbf{E}_x$ the corresponding expectation.
It is well-known that the distribution of $Y$ can be specified by $\Phi$ as follows: For $\lambda\geq 0$,
$$
\mathbf{E}_x[\exp(-\lambda Y_t)]=\exp[-x v_t(\lambda)],
$$
where $v_t(\lambda)$ is the unique locally bounded nonnegative solution of
$$
v_t(\lambda)=-\int_0^t \Phi (v_s(\lambda)) ds +\lambda.
$$

The \emph{excursion representation} of CB processes is important in this paper.
Take a CB process $Y$ with the branching mechanism $\Phi$, we can define an excursion measure $N$ and reconstruct $Y$ from excursions.
Let $\mathbb{D}([0,\infty),\mathbb{R}_+)$ be the standard Skorohod's space.
Let $\mathbb{D}_0([0,\infty),\mathbb{R}_+)$ be the subspace of $\mathbb{D}([0,\infty),\mathbb{R}_+)$, such that all paths in
$\mathbb{D}_0([0,\infty),\mathbb{R}_+)$ start from $0$ and stop upon hitting $0$.
Under Condition (\ref{assumptionweak}), we may define a $\sigma$-finite measure $N$ on $\mathbb{D}_0([0,\infty),\mathbb{R}_+)$ such that:\\
\noindent 1. $N(\{\mathbf{0}\})=0$, where $\mathbf{0}$ denotes the trivial path in $\mathbb{D}([0,\infty),\mathbb{R}_+)$, that is,
$\mathbf{0}_t=0$ for any $t$.\\
2. Let $Z$ be a Poisson random measure on $\mathbb{D}_0([0,\infty),\mathbb{R}_+)$ with intensity $xN$.
Define the process $(e_t,t\geq 0)$ by $e_0= x$ and
$$
e_t=\int_{\mathbb{D}_0([0,\infty),\mathbb{R}_+)}\omega_t Z(d\omega), \quad t>0.
$$
Then $e$ is a CB process with the branching mechanism $\Phi$.

\subsection{Height processes and the branching property of L\'{e}vy trees}\label{ss:bp}

The \emph{height process} $H=(H_t,t\geq 0)$ is introduced by Le Gall and Le Jan \cite{LL98} and further developed by Duquesne and Le Gall \cite{DL02},
to code the complete genealogy of L\'{e}vy trees. It is obtained as a functional of the spectrally positive L\'{e}vy process $X$ with
the Laplace exponent $-\Phi$.
Intuitively, for every $t \geq 0$, $H_t$ ``measures" in a local time sense the size of
the set $\{s \leq t: X_{s-}=\inf_{r\in[s,t]}X_r\}$.
Condition (\ref{assumptionstrong}) holds if and only if
$H$ has a continuous modification. From now on, we only consider this modification.
For any $a\geq 0$, the local time $L^a=(L^a_t,t\geq 0)$ of $H$ at height $a$ can be defined, which is continuous and increasing. Intuitively, the measure induced by $L^a$ is distributed ``uniformly" on all ``particles" of the L\'{e}vy tree at height $a$. Note that all processes introduced so far are defined under the underlying probability $\mathbf{P}$, so that all these processes correspond to a Poisson collection of L\'{e}vy trees.
This Poisson collection has infinite but $\sigma$-finite intensity, which in turn corresponds to a CB process with infinite initial value.

In order to talk about a single L\'{e}vy tree, a certain excursion measure $\mathbf{N}$ needs to be introduced. Recall the spectrally positive L\'{e}vy process $X=(X_t,t\geq 0)$ with the Laplace exponent $-\Phi$, and its infimum process $I=(I_t,t\geq 0)$ defined by
$I_t=\inf_{s\leq t}X_s$.
When Condition (\ref{assumptionweak}) holds, the point 0 is regular and instantaneous for the strong Markov process $X-I$. We denote by $\mathbf{N}$ the corresponding excursion
measure, and by $\zeta$ the duration of the excursion. We also denote by $X$ the canonical process under $\mathbf{N}$.
Note that normally we need to specify the normalization of $\mathbf{N}$, but for our purposes in this paper this normalization always cancels out.

In general $H$ is not Markov under $\mathbf{P}$, but $H_t$ only depends on the values of $X-I$, on the
excursion interval of $X-I$ away from 0 that straddles $t$. Also it can be checked that a.s. for all $t$, $H_t>0$ if and only if $X_t-I_t>0$.
So under $\mathbf{N}$ we may define $H$ as a functional of $X$ (recall that $X$ is the canonical process under $\mathbf{N}$). Consequently we may also define $L^a=(L^a_t,0\leq t\leq \zeta)$ of $H$ at any height $a\geq 0$, under the excursion measure $\mathbf{N}$. Note that it is then standard to define the L\'{e}vy tree with the branching mechanism $\Phi$ as a random metric space $\mathcal{T}(\Phi)$ from the height process $H$ and in this paper we will just regard $H$ as the L\'{e}vy tree.

The \emph{branching property} of L\'{e}vy trees is crucial for us in this paper. For any $b>0$, define the conditional probability $\mathbf{N}^{(b)}$
as the distribution of the canonical process $X$ conditioned on having height greater than $b$, that is,
$$\mathbf{N}^{(b)}[\cdot]=\mathbf{N}[\cdot|\sup H>b].$$
Then intuitively the branching property says that under $\mathbf{N}^{(b)}$ and conditional on all information below height $b$, all the subtrees above height $b$ are just i.i.d. copies of the complete L\'{e}vy tree under $\mathbf{N}$, and the roots of all these subtrees distribute as a Poisson random measure  with intensity the measure induced by $L^b=(L^b_t,0\leq t\leq \zeta)$.
It is well-known that $\mathbf{N}[(L^b_\zeta,b\geq 0)\in \cdot]=N[\cdot]$.
So consequently from the excursion representation of CB processes, we see that under $\mathbf{N}^{(b)}$ and conditional on all information below height $b$, the real-valued process $(L^a_\zeta,b\leq a < \infty)$ distributes as a CB process with initial value $L^b_\zeta$.
For a rigorous presentation of this branching property, refer to Proposition 3.1 in \cite{D08} or Corollary 3.2 in \cite{DL05}.

\subsection{Continuum immortal trees}\label{ss:CKT}

Recall (\ref{Levy}) and (\ref{subordinator}).
Let $H$ be the height process associated with the L\'{e}vy process $X$ with the Laplace exponent $-\Phi$ and let $(H',X')$ be a copy of $(H,X)$.
Let $I=(I_t,t\geq 0)$ and $I'=(I'_t,t\geq 0)$ be the infimum processes of $X$ and $X'$ respectively.
Let $(U,V)$ be a bivariate subordinator with the Laplace exponent $\Phi(p,q)$.
Let $U^{-1}=(U^{-1}_t,t\geq 0)$ and $V^{-1}=(V^{-1}_t,t\geq 0)$ be the right-continuous inverses of $U$ and $V$ respectively.
Assume that $(X,H)$, $(X',H')$, and $(U,V)$ are independent.
We define $\overleftarrow{H}$ and $\overrightarrow{H}$ by
$$
\overleftarrow{H}_t=H_t+U^{-1}_{-I_t} \quad \text{and}\quad
\overrightarrow{H}_t=H'_t+V^{-1}_{-I'_t}, \quad t\geq 0.
$$
The processes $\overleftarrow{H}$ and $\overrightarrow{H}$ are called respectively left and right height processes of the continuum immortal tree with the branching mechanism $\Phi$.
Then it is natural to define the continuum immortal tree as a random metric space $\mathcal{T}^*(\Phi)$ from the height processes $\overleftarrow{H}$ and $\overrightarrow{H}$.
For details refer to Page 103 in \cite{D08}. In this paper we will just regard the height processes $\overleftarrow{H}$ and $\overrightarrow{H}$ as the continuum immortal tree.

Introduce the last time under level $b$ for the left
and the right height processes:
$$
\overleftarrow{\sigma}_b=\sup\{s\geq 0: \overleftarrow{H}_s\leq b\}\quad \text{and} \quad \overrightarrow{\sigma}_b=\sup\{s\geq 0: \overrightarrow{H}_s\leq b\}.
$$
Now let us recall Lemma 3.2 in \cite{D08}, which relates the distribution of a L\'{e}vy tree and that of the corresponding continuum immortal tree.

\begin{lem} \label{Bismut}
For any nonnegative measurable functions $F$ and $G$, and any $b>0$,
$$
\mathbf{N}\l[\int_0^\zeta \text{d} L^b_s F\l(H_{\cdot\wedge s}\r)G\l(H_{(\zeta-\cdot)\wedge (\zeta-s)}\r)\r]=e^{-\alpha b}\mathbf{E}\l[F(\overleftarrow{H}_{\cdot\wedge \overleftarrow{\sigma}_b})G(\overrightarrow{H}_{\cdot\wedge \overrightarrow{\sigma}_b})\r].
$$
\end{lem}

Note that taking $F=G\equiv 1$ in Lemma \ref{Bismut} gives
\beqlb \label{expectation}
\mathbf{N}[L^b_\zeta]=e^{-\alpha b}, \quad b>0.
\eeqlb

\section{Local convergence of critical L\'{e}vy trees and CB processes} \label{s:C}

In this section first we study the local convergence of conditioned critical L\'{e}vy trees.
Then we derive a very general ratio limit property on certain functionals of L\'{e}vy trees.
Finally we treat the local convergence of conditioned critical CB processes.

\subsection{Local convergence of critical L\'{e}vy trees}\label{ss:L}
For $\omega=(\omega_t,t\geq 0)\in\mathbb{C}([0,\infty),\mathbb{R}_+)$, define $\zeta(\omega)=\inf\{t>0:\omega_t=0\}$.
Denote by $\mathbb{C}_0([0,\infty),\mathbb{R}_+)$ the subspace of all excursions in $\mathbb{C}([0,\infty),\mathbb{R}_+)$, that is, $\omega\in\mathbb{C}_0([0,\infty),\mathbb{R}_+)$ if and only if $\omega\in\mathbb{C}([0,\infty),\mathbb{R}_+)$, $\omega_t>0$ when $t\in(0,\zeta)$, and $\omega_t=0$ otherwise.

Let $A$ be a nonnegative measurable function defined on $\mathbb{C}_0([0,\infty),\mathbb{R}_+)$.
For an excursion $\omega\in\mathbb{C}_0([0,\infty),\mathbb{R}_+)$, write $\omega_{(b)}=(\omega_{(b),i},i\in \mathcal{I}_{(b)})$ for the collection of all sub-excursions above height $b$.
We introduce the following monotonicity property of $A$:
\beqlb \label{Lmonotonicity}
A(\omega_{(b),i})\leq A(\omega), \quad \text{for any } \omega\in\mathbb{C}_0([0,\infty),\mathbb{R}_+),\ b>0, \text{ and }i\in \mathcal{I}_{(b)}.
\eeqlb
Suppose that $\omega\in\mathbb{C}_0([0,\infty),\mathbb{R}_+)$ is the height process of a real tree, then the monotonicity property (\ref{Lmonotonicity}) says that for any $b>0$ the value of $A$ on the whole tree is not less than that on any subtree above height $b$.

Define $v_r=\mathbf{N}[A (H)>r]\in [0,\infty]$ and $\mathbf{N}_r[\cdot]=\mathbf{N}[\cdot|A (H)>r]$ when $v_r\in (0,\infty)$.
The following theorem asserts that if the monotonicity property (\ref{Lmonotonicity}) holds for $A$, then
under the conditional probability $\mathbf{N}_r$, the L\'{e}vy tree $\mathcal{T}(\Phi)$ converges locally to the continuum immortal tree $\mathcal{T}^*(\Phi)$, see Remark \ref{tree}.

\begin{theo} \label{L}
Assume that $\Phi$ is critical and $v_r\in (0,\infty)$ for large enough $r$.
If the function $A$ satisfies the monotonicity property (\ref{Lmonotonicity}),
then as $r\rightarrow\infty$,
$$
(H_{t\wedge \zeta},H_{(\zeta-t)_+};t\geq 0)
\text{ under }\mathbf{N}_r\longrightarrow
(\overleftarrow{H}_t,\overrightarrow{H}_t;t\geq 0)
$$
weakly in $\mathbb{C}([0,\infty),\mathbb{R}^2)$.
\end{theo}

\proof First we follow the beginning of the proof of Theorem 1.3 in \cite{D08} and for the reader's convenience we copy that part here.
Let $b > 0$. For any $\omega$ in $\mathbb{C}_0([0,\infty),\mathbb{R}_+)$ we introduce $\tau_b(\omega)=\inf \{s\geq 0: \omega(s)=b\}$.
To simplify notations we set
$$
\widehat{H}_\cdot=H_{(\zeta-\cdot)_+},\quad \tau_b=\tau_b(H), \quad\widehat{\tau}_b=\tau_b(\widehat{H}), \quad
\overleftarrow{\tau}_b=\tau_b(\overleftarrow{H}),
\quad\text{and}\quad \overrightarrow{\tau}_b=\tau_b(\overrightarrow{H}).
$$
We only have to prove the following convergence for any bounded measurable function $F$,
\beqlb \label{remark}
\lim_{r\rightarrow\infty} \mathbf{N}_r\l[F(H_{\cdot\wedge \tau_b},\widehat{H}_{\cdot\wedge \widehat{\tau}_b})\r]
=\mathbf{E}\l[F(\overleftarrow{H}_{\cdot\wedge \overleftarrow{\tau}_b},\overrightarrow{H}_{\cdot\wedge \overrightarrow{\tau}_b})\r],\quad b>0,
\eeqlb
since it implies for any $t > 0$,
$$
\lim_{b\rightarrow\infty}\lim_{r\rightarrow\infty} \mathbf{N}_r[\tau_b\wedge \widehat{\tau}_b\leq t ]
=\lim_{b\rightarrow\infty}\mathbf{P}\l[\overleftarrow{\tau}_b\wedge \overrightarrow{\tau}_b\leq t\r]=0.
$$
We may deduce from a standard approximation result of $L^b$ (see (28) in \cite{D08}) that $\mathbf{N}$ a.e. the topological support of $dL^b_\cdot$ is included in
$[\tau_b,\zeta-\widehat{\tau}_b]$. Thus, by Lemma \ref{Bismut},
$$
\mathbf{N}\l[L^b_\zeta F(H_{\cdot\wedge \tau_b},\widehat{H}_{\cdot\wedge \widehat{\tau}_b})\r]=\mathbf{E}\l[F(\overleftarrow{H}_{\cdot\wedge \overleftarrow{\tau}_b},\overrightarrow{H}_{\cdot\wedge \overrightarrow{\tau}_b})\r],\quad b>0.
$$
So it suffices to show that
\beqlb \label{enough}
\lim_{r\rightarrow\infty} \mathbf{N}_r\l[F(H_{\cdot\wedge \tau_b},\widehat{H}_{\cdot\wedge \widehat{\tau}_b})\r]
=\mathbf{N}\l[L^b_\zeta F(H_{\cdot\wedge \tau_b},\widehat{H}_{\cdot\wedge \widehat{\tau}_b})\r],\quad b>0.
\eeqlb

To prove (\ref{enough}), first recall that if the value of $A$ on a subtree above height $b$ is greater than $r$, than the value of $A$ on the whole tree is greater than $r$, by the monotonicity property (\ref{Lmonotonicity}). Then recall the branching property from Section \ref{ss:bp} and
note that under $\mathbf{N}_{(b)}$ and conditional on all information below height $b$ the probability that the value of $A$ on at least one subtree above height $b$ is greater than $r$ is
$$1-\exp\l(-L^b_\zeta v_r\r).$$
So the monotonicity property and the branching property imply that
$$
\mathbf{N}^{(b)}\l[\mathbf{1}\{A (H)>r\}  F(H_{\cdot\wedge \tau_b},\widehat{H}_{\cdot\wedge \widehat{\tau}_b})  \r]\geq
\mathbf{N}^{(b)}\l[\l(1-e^{-L^b_\zeta v_r}\r)  F(H_{\cdot\wedge \tau_b},\widehat{H}_{\cdot\wedge \widehat{\tau}_b})  \r],
$$
which further implies that
$$
\mathbf{N}\l[\mathbf{1}\{A (H)>r\}  F(H_{\cdot\wedge \tau_b},\widehat{H}_{\cdot\wedge \widehat{\tau}_b})  \r]\geq
\mathbf{N}\l[\l(1-e^{-L^b_\zeta v_r}\r)  F(H_{\cdot\wedge \tau_b},\widehat{H}_{\cdot\wedge \widehat{\tau}_b})  \r].
$$
Thus we see that as $r\rightarrow\infty$,
\beqnn
\mathbf{N}_r\l[F(H_{\cdot\wedge \tau_b},\widehat{H}_{\cdot\wedge \widehat{\tau}_b})\r]
\ar\geq\ar\frac{1}{v_r}\mathbf{N}\l[\l(1-e^{-L^b_\zeta v_r}\r)  F(H_{\cdot\wedge \tau_b},\widehat{H}_{\cdot\wedge \widehat{\tau}_b})  \r] \cr
\ar\rightarrow\ar \mathbf{N}\l[L^b_\zeta F(H_{\cdot\wedge \tau_b},\widehat{H}_{\cdot\wedge \widehat{\tau}_b})\r],\cr
\eeqnn
where the convergence follows from the monotone convergence.

From the above paragraph we get the inequality that
$$
\liminf_{r\rightarrow\infty} \mathbf{N}_r\l[F(H_{\cdot\wedge \tau_b},\widehat{H}_{\cdot\wedge \widehat{\tau}_b})\r]
\geq \mathbf{N}\l[L^b_\zeta F(H_{\cdot\wedge \tau_b},\widehat{H}_{\cdot\wedge \widehat{\tau}_b})\r].
$$
Clearly we may assume that $0\leq F\leq 1$, then apply the above inequality to $1-F$ to get
$$
\liminf_{r\rightarrow\infty} \mathbf{N}_r\l[1-F(H_{\cdot\wedge \tau_b},\hat{H}_{\cdot\wedge \hat{\tau}_b})\r]
\geq \mathbf{N}\l[L^b_\zeta -L^b_\zeta F(H_{\cdot\wedge \tau_b},\widehat{H}_{\cdot\wedge \widehat{\tau}_b})\r],
$$
which implies that
$$
\limsup_{r\rightarrow\infty} \mathbf{N}_r\l[F(H_{\cdot\wedge \tau_b},\hat{H}_{\cdot\wedge \hat{\tau}_b})\r]
\leq \mathbf{N}\l[L^b_\zeta F(H_{\cdot\wedge \tau_b},\widehat{H}_{\cdot\wedge \widehat{\tau}_b})\r],
$$
since $\mathbf{N}[L^b_\zeta]=1$ by (\ref{expectation}). Finally we have proved (\ref{enough}).
\qed

\begin{rem} \label{tree}
Recall that the weak convergence on $\mathbb{C}([0,\infty),\mathbb{R}^2)$ is defined with respect to the topology of uniform convergence on compact subsets of $[0,\infty)$. This corresponds to the local convergence of random real trees that we consider here. More specifically, recall that in the proof of Theorem \ref{L} we have proved (\ref{remark}). By combining (\ref{remark}) with Lemma 2.3 in \cite{DL05}, we see that our Theorem \ref{L} says that under the conditional probability
$\mathbf{N}_r$, the subtree of the L\'{e}vy tree $\mathcal{T}(\Phi)$ below height $b$ converges to the subtree of the continuum immortal tree $\mathcal{T}^*(\Phi)$ below height $b$, with respect to the Gromov-Hausdorff distance on the space of all equivalence
classes of rooted compact real trees.
\end{rem}

Next we will apply Theorem \ref{L} to three specific conditionings, which are the conditioning of large width, the conditioning of large total mass, and the conditioning of large maximal degree.
We first introduce the conditioning of large width.
Under $\mathbf{N}$, define the \emph{width} of the L\'{e}vy tree $H$ by $W(H)=\sup_{b\geq 0}L^b_\zeta$.
Consider $\mathbf{N}[\cdot|W(H)>r]$ when $\mathbf{N}[W (H)>r]\in (0,\infty)$, this is the conditioning of large width.
Then the conditioning of large total mass. Under $\mathbf{N}$, define the \emph{total mass} of the L\'{e}vy tree $H$ by $\sigma(H)=\int_0^\infty L^b_\zeta db$.
Consider $\mathbf{N}[\cdot|\sigma(H)>r]$ when $\mathbf{N}[\sigma(H)>r]\in (0,\infty)$, this is the conditioning of large total mass.

Finally we introduce the conditioning of large maximal degree. Recall from Section \ref{ss:bp} that $\mathbf{N}$ is the excursion
measure of the strong Markov process $X-I$ at zero.
Also recall that we write $X$ for the canonical process under $\mathbf{N}$, which is rcll.
Finally recall from Theorem 4.6 of \cite{DL05} that L\'{e}vy trees have two types of nodes (i.e., branching points), binary nodes
(i.e., vertices of degree 3) and infinite nodes (i.e., vertices of infinite degree). Infinite nodes correspond to the
jumps of the canonical process $X$ under $\mathbf{N}$, and the sizes of these jumps correspond to the masses of those infinite nodes.
We call the mass of a node its \emph{degree}.
Then define the \emph{maximal degree} of the L\'{e}vy tree $H$ by $M(H)=\sup_{0\leq s\leq \zeta}\Delta X_s$. Note that under $\mathbf{N}$ we can write $\max_{0\leq s\leq \zeta}\Delta X_s$ as a functional of $H$, since jumps of $X$ correspond to jumps of $L_\zeta=(L^b_\zeta,b\geq 0)$, which are functionals of $H$.
Consider $\mathbf{N}[\cdot|M(H)>r]$ when $\mathbf{N}[M (H)>r]\in (0,\infty)$, this is the conditioning of large maximal degree.

Since the monotonicity property (\ref{Lmonotonicity}) is trivial to check, we then only have to check that $v_r\in (0,\infty)$ for large enough $r$. In the following lemma, we only assume (\ref{assumptionweak}). Note that to define $W$, $\sigma$, and $M$, we only need the real-valued process $(L^a_\zeta,a\geq 0)$, which when (\ref{assumptionweak}) holds can be defined by the excursion representation of CB processes, without the introduction of the height process $H$ and its local times. Also note that the functionals $W$, $\sigma$, or $M$ can be similarly defined for CB processes. We write $\mathbf{P}_x$ for probabilities of CB processes with initial value $x$.

\begin{lem} \label{positive} For any $x>0$, if $\alpha\geq 0$, then $\mathbf{P}_x[W >r]>0$ and $\mathbf{N}[W>r]\in (0,\infty)$ for any $r\in (0,\infty)$, and $\mathbf{P}_x[\sigma >r]>0$ and $\mathbf{N}[\sigma>r]\in (0,\infty)$ for any $r\in (0,\infty)$. Again for any $x>0$, if $\alpha\geq 0$ and $\pi$ has unbounded support, then $\mathbf{P}_x[M >r]>0$ and $\mathbf{N}[M>r]\in (0,\infty)$ for any $r\in (0,\infty)$.
\end{lem}

\proof
For the width, first we argue that for any $x\in (0,\infty)$ and $r\in (0,\infty)$,
$$
\mathbf{P}_x[W >r]\leq \frac{x}{r}.
$$
When $\alpha\geq 0$, for the CB process $Y$ we may define $Y_\infty=0$ and regard $Y$ as a supermartingale over the time interval $[0,\infty]$.
Then by optional sampling, it is easy to get the above inequality.
Now by the excursion representation of CB processes and the above inequality, we get
$$
1-\exp\l(-(r/2)\mathbf{N}[W>r]\r)\leq \mathbf{P}_{r/2}[W>r]\leq 1/2,
$$
which implies that $\mathbf{N}[W>r]<\infty$ for any $r>0$.
Next by Corollary 12.9 in \cite{K14} and the fact that scale functions are strictly increasing, we see that $\mathbf{P}_x[W>r]>0$ for any $x>0$ and $r>0$.
Finally by the branching property of L\'{e}vy trees, for any $r>0$ and $b>0$ we have
$$
\mathbf{N}[W>r]\geq\mathbf{N}\l[\mathbf{P}_{L^b_\zeta}[W>r]\r]>0,
$$
since $\mathbf{N}[L^b_\zeta>0]>0$ for any $b>0$.

For the total mass, first denote by $\Phi^{-1}$ the inverse function of $\Phi$,
then recall that for $\lambda>0$, $\mathbf{N}[1-\exp(-\lambda\sigma)]=\Phi^{-1}(\lambda)$ (see e.g., the beginning of Section 3.2.2 in \cite{DL02}), which implies that
$\mathbf{N}[\sigma>r]<\infty$ for any $r>0$. Also clearly $\mathbf{N}[\sigma>r]>0$ for some $r>0$, then by the excursion representation of CB processes, we see that $\mathbf{P}_x[\sigma>r]>0$ for any $r>0$. Finally by the branching property of $\mathbf{N}$, for any $r>0$ and $b>0$ we have
$$
\mathbf{N}[\sigma>r]\geq\mathbf{N}\l[\mathbf{P}_{L^b_\zeta}[\sigma>r]\r]>0.
$$

For the maximal degree, trivially $\mathbf{N}[M>r]\leq\mathbf{N}[W>r]<\infty$ for any $r\in (0,\infty)$.
Then by the one-to-one correspondence between distributions of CB processes and branching mechanisms, we know that for any $r\in (0,\infty)$,
$\mathbf{P}_x[M  >r]>0$.
Finally by the branching property of L\'{e}vy trees, for any $r\in (0,\infty)$ and $b\in (0,\infty)$ we have
$$
\mathbf{N}[M>r]\geq\mathbf{N}\l[\mathbf{P}_{L^b_\zeta}[M>r]\r]>0.
$$
\qed

Now Theorem \ref{L} and Lemma \ref{positive} immediately imply the local convergence of critical L\'{e}vy trees, under any of the three conditionings we introduced above.

\begin{cor} \label{L3}
Assume that $\Phi$ is critical.
Then for $A=W$ and $A=\sigma$ respectively, as $r\rightarrow\infty$,
$$
(H_{t\wedge \zeta},H_{(\zeta-t)_+};t\geq 0)
\text{ under }\mathbf{N}[\cdot|A(H)>r]\longrightarrow
(\overleftarrow{H}_t,\overrightarrow{H}_t;t\geq 0)
$$
weakly in $C([0,\infty),\mathbb{R}^2)$. The above local convergence also holds for $A=M$ under the additional assumption that
$\pi$ has unbounded support.
\end{cor}

\subsection{A general ratio limit property}\label{ss:CRL}
Denote by $\mathbb{C}_0^\infty([0,\infty),\mathbb{R}_+)$ the product space of countably infinite copies of $\mathbb{C}_0([0,\infty),\mathbb{R}_+)$.
Let $A$ be a nonnegative measurable function defined on $\mathbb{C}_0^\infty([0,\infty),\mathbb{R}_+)$, which is invariant under permutation.
For any $\omega\in\mathbb{C}_0([0,\infty),\mathbb{R}_+)$, we also write
$A(\omega)=A(\omega^\infty)$, where $\omega^\infty=(\omega,\mathbf{0},\mathbf{0},\ldots)\in\mathbb{C}_0^\infty([0,\infty),\mathbb{R}_+)$.
We introduce the following monotonicity property of $A$:
\beqlb \label{Cmonotonicity}
A(\omega_{(b),i})\leq A(\omega_{(b)})\leq A(\omega), \quad \text{for any } \omega\in\mathbb{C}_0([0,\infty),\mathbb{R}_+),\ b>0, \text{ and }i\in \mathcal{I}_{(b)}.
\eeqlb

Let $N^{(x)}$ be a Poisson random measure on $\mathbb{C}_0([0,\infty),\mathbb{R}_+)$ with intensity $x\mathbf{N}[H\in\cdot]$. Write $\omega^{(x)}=(\omega^{(x),i},i\in \mathcal{I}^{(x)})$ for the collection of all excursions in $N^{(x)}$.
Also write $\mathbf{N}[A>r]$ for $\mathbf{N}[A(H)>r]$, and $\mathbf{P}^{(x)}[A>r]$ for $\mathbf{P}[A(\omega^{(x)})>r]$.
Define $r_b(\omega)=((r_b(\omega))_s,s\geq 0)$ by $(r_b(\omega))_s=\omega_{s_b}$, where
$$
s_b=\inf\{t\geq 0:\int_0^t da\mathbf{1}\{H_a\leq b\}>s\}.
$$
So if $\omega\in\mathbb{C}_0([0,\infty),\mathbb{R}_+)$ is the height process of a real tree, then $r_b(\omega)$ is just the height process of the corresponding subtree below height $b$.
Write
$$\mathbb{C}_0^{(b)}([0,\infty),\mathbb{R}_+)=\{r_b(\omega):\omega\in\mathbb{C}_0([0,\infty),\mathbb{R}_+)\}.$$
We also introduce an additivity property of $A$: For any fixed $b>0$ and $\omega' \in\mathbb{C}_0^{(b)}([0,\infty),\mathbb{R}_+)$,
\beqlb \label{Cadditivity}
A(\omega)=A(\omega_{(b)})+ B(r_b(\omega)), \quad \text{for large enough }A(\omega)\text{ with }r_b(\omega)=\omega',
\eeqlb
where $B$ is a nonnegative measurable function defined on $\mathbb{C}_0([0,\infty),\mathbb{R}_+)$.

\begin{theo} \label{CRL}
Assume that $\Phi$ is critical and $\mathbf{N}[A>r]\in(0,\infty)$ for large enough $r$.
If $A$ satisfies the monotonicity property (\ref{Cmonotonicity}),
then for any $x>0$,
$$
\lim_{r\rightarrow\infty} \frac{\mathbf{P}^{(x)}[A>r]}{\mathbf{N}[A>r]}=x.
$$
If $A$ satisfies the additivity property (\ref{Cadditivity}), and for some $b>0$ and any $\omega' \in\mathbb{C}_0^{(b)}([0,\infty),\mathbb{R}_+)$,
$B(r_b(\omega))>0$ for large enough $A(\omega)$ with $r_b(\omega)=\omega'$,
then for any $x>0$ and $r'>0$,
$$
\lim_{r\rightarrow\infty} \frac{\mathbf{P}^{(x)}[A>r-r']}{\mathbf{P}^{(x)}[A>r]}=1.
$$
\end{theo}

\proof First as in the proof of Theorem \ref{L}, for any $x>0$,
$$
\liminf_{r\rightarrow\infty} \frac{\mathbf{P}^{(x)}[A>r]}{\mathbf{N}[A>r]}\geq \liminf_{r\rightarrow\infty} \frac{1-\exp(-x\mathbf{N}[A>r])}{\mathbf{N}[A>r]}= x.
$$

Next we show that there exists some unbounded $K\subset\mR_+$, such that for any $x\in K$,
\beqlb \label{=}
\lim_{r\rightarrow\infty} \frac{\mathbf{P}^{(x)}[A>r]}{\mathbf{N}[A>r]}= x.
\eeqlb
To prove this, just note that as in the proof of Theorem \ref{L}, we have
$$
\lim_{r\rightarrow\infty}\frac{\mathbf{N}\l[ \mathbf{P}^{(L^b_\zeta)}[A>r] F(H_{\cdot\wedge \tau_b},\widehat{H}_{\cdot\wedge \widehat{\tau}_b})  \r]}{\mathbf{N}[A>r]}=\mathbf{N}\l[L^b_\zeta F(H_{\cdot\wedge \tau_b},\widehat{H}_{\cdot\wedge \widehat{\tau}_b})\r],
$$
which implies that for a.e. $x$ with respect to dist$(L^b_\zeta)$,
$$
\lim_{r\rightarrow\infty}\frac{\mathbf{P}^{(x)}[A>r]}{\mathbf{N}[A>r]}= x.
$$
Recall from Lemma \ref{positive} that $W=\sup_{b\in\mR_+}L^b_\zeta$ has unbounded support, then it is immediate that the desired $K$ exists.

Finally assume that for some $x>0$,
$$
\limsup_{r\rightarrow\infty} \frac{\mathbf{P}^{(x)}[A>r]}{\mathbf{N}[A>r]}>x.
$$
Clearly we can pick some $x'\in K$ such that $x'>x$, so that (\ref{=}) holds for $x'$.
However, we also have
\beqnn
\mathbf{P}^{(x')}[A>r]\ar\geq\ar 1-(1-\mathbf{P}^{(x)}[A>r])\exp(-(x'-x)\mathbf{N}[A>r])\cr
\ar=\ar 1-\exp(-(x'-x)\mathbf{N}[A>r])+\mathbf{P}^{(x)}[A>r]\exp(-(x'-x)\mathbf{N}[A>r]),\cr
\eeqnn
which implies that
$$
\limsup_{r\rightarrow\infty} \frac{\mathbf{P}^{(x')}[A>r]}{\mathbf{N}[A>r]}
\geq (x'-x)+\limsup_{n\rightarrow\infty}\frac{\mathbf{P}^{(x)}[A>r]}{\mathbf{N}[A>r]}>x',
$$
a contradiction to (\ref{=}) for $x'$.

For the second statement, by the argument in the proof of Theorem \ref{GWtail}, the assumptions, and a disintegration theorem (see e.g., Theorem 6.4 in \cite{K02}), we have
$$
\lim_{n\rightarrow\infty}\frac{\mathbf{N}\l[\mathbf{P}^{(L^b_\zeta)}[A>r-B(r_b(H))]F(H_{\cdot\wedge \tau_b},\widehat{H}_{\cdot\wedge \widehat{\tau}_b})\r]}{\mathbf{N}[A>r]}=\mathbf{N}\l[L^b_\zeta F(H_{\cdot\wedge \tau_b},\widehat{H}_{\cdot\wedge \widehat{\tau}_b})\r]
$$
which means that there exist some $x>0$ and $r'>0$ such that
$$
\lim_{r\rightarrow\infty} \frac{\mathbf{P}^{(x)}[A>r-r']}{\mathbf{N}[A>r]}=x.
$$
Combined with the first statement, we get for this particular $r'>0$,
$$
\lim_{r\rightarrow\infty} \frac{\mathbf{N}[A>r-r']}{\mathbf{N}[A>r]}=1,
$$
which implies that for any $r'>0$,
$$
\lim_{r\rightarrow\infty} \frac{\mathbf{N}[A>r-r']}{\mathbf{N}[A>r]}=1.
$$
Finally use the first statement again to finish the proof of the second statement.
\qed

Next we give some applications of Theorem \ref{CRL} (or rather a variant of Theorem \ref{CRL}, see the explanation after the following corollary).
The first application is about the width of CB processes and the scale functions of L\'{e}vy processes. For the definition of the scale function $W=(W(r),r\geq0)$, see e.g., Section 8.3 in \cite{K14}.
The second application is about the total mass of CB processes.
The following corollary is immediate from our Lemma \ref{positive} and Theorem \ref{CRL} (or more precisely the variant of Theorem \ref{CRL} we just mentioned), and Corollary 12.9 in \cite{K14}.
It might be interesting to note that the convergence of the scale function $W(r)$ below also has an intuitive meaning for first passage times of L\'{e}vy processes, see (8.11) in \cite{K14}.

\begin{cor} \label{applications}
Assume that $\Phi$ is critical and satisfies (\ref{assumptionweak}),
then for any $x>0$,
$$
\lim_{r\rightarrow\infty} \frac{\mathbf{P}_x[W>r]}{N[W>r]}=x.
$$
Expressed in terms of the scale function $W=(W(r),r\geq 0)$, the above convergence means that for any $x>0$,
$$
\lim_{r\rightarrow\infty} \frac{W(r)-W(r-x)}{W(r)-W(r-1)}=x.
$$
Also for any $x>0$,
$$
\lim_{r\rightarrow\infty} \frac{\mathbf{P}_x[\sigma>r]}{N[\sigma>r]}=x,
$$
and for any $x,r'>0$,
$$
\lim_{r\rightarrow\infty} \frac{\mathbf{P}_x[\sigma>r-r']}{\mathbf{P}_x[\sigma>r]}=1.
$$
\end{cor}

Note that Theorem \ref{CRL} is about height processes of L\'{e}vy trees, so we have to assume Condition (\ref{assumptionstrong}) to get continuous height processes. However the proof does not rely on any specific property of height processes.
In fact if we are only interested in the real-valued process $(L^a_\zeta,a\geq 0)$, not the height process and its local times, then as long as the branching mechanism satisfies the weaker Condition (\ref{assumptionweak}), we can make the same proof work by using the excursion representation of CB processes.
Recall that $W$, $\sigma$, and $M$ are all functionals of the real-valued process $(L^a_\zeta,a\geq 0)$, so in Corollary \ref{applications} we only need to assume Condition (\ref{assumptionweak}). Also note that in Corollary \ref{applications} we use $\mathbf{P}_x$ and $N$ instead of $\mathbf{P}^{(x)}$ and $\mathbf{N}$, respectively.

We do not need to study the ratio limit property of the maximal degree of L\'{e}vy forests since it is trivial.
Assume that for any $\omega^\infty=(\omega_1,\omega_2,\ldots)\in\mathbb{C}_0^\infty([0,\infty),\mathbb{R}_+)$,
the functional $A$ has the property that
\beqlb \label{max}
A(\omega^\infty)=\sup_{i\in\mN} A(\omega_i).
\eeqlb
Then clearly
$$
\mathbf{P}^{(x)}[A>r]=1-\exp(-x\mathbf{N}[A>r]),
$$
which implies that
$$
\lim_{r\rightarrow\infty} \frac{\mathbf{P}^{(x)}[A>r]}{\mathbf{N}[A>r]}=x.
$$
The maximal degree and the height of L\'{e}vy forests are two examples of functionals satisfying (\ref{max}). Clearly for the maximal jump and the height of CB processes we have a similar situation.

\subsection{Local convergence of critical CB processes}\label{ss:LCCB}

In this sub-section we treat the local convergence of conditioned critical CB processes.
To avoid repetitions of several arguments in the previous sub-sections, here we only set up the framework and state the result.

Let $A$ be a nonnegative measurable function defined on $\mathbb{D}([0,\infty),\mathbb{R}_+)$, the standard Skorohod's space.
For any $\omega=(\omega_t,t\geq 0)\in\mathbb{D}([0,\infty),\mathbb{R}_+)$, write $\omega_{(b)}=(\omega_t,t\geq b)$ for the sub-path after time $b$.
Note that this $\omega_{(b)}$ defined here corresponds to the $\omega_{(b)}$ defined in Section \ref{ss:L}, so we slightly abuse the notation.
We introduce the following monotonicity property of $A$:
\beqlb \label{CBmonotonicity}
A(\omega_{(b)})\leq A(\omega), \quad \text{for any } \omega\in\mathbb{D}([0,\infty),\mathbb{R}_+)\text{ and } b>0.
\eeqlb

Let $Y$ be a CB process with the branching mechanism $\Phi$ and
write $(\mathcal{F}_b,b\geq 0)$ for the filtration induced by $Y$.
Use $\mathbf{P}_x[Y\in \cdot]$ to denote the distribution of $Y$ with $Y_0=x$.
The following theorem can be proved by adapting the proofs of Theorem \ref{L} and Theorem \ref{CRL} (note that we need Condition (\ref{assumptionweak}) to use the excursion representation of CB processes).
It asserts that if the monotonicity property (\ref{CBmonotonicity}) holds for $A$, then
under the conditional probability $\mathbf{P}_x[\cdot|A (Y)>r]$, the CB process $Y$ converges locally to a CB process with immigration (CBI process), such that the branching mechanism of this CBI process is still $\Phi$ and the immigration mechanism is $\Phi'$, the derivative of $\Phi$,
see Remark \ref{CBI}.

\begin{theo} \label{LCCB}
Assume that $\Phi$ is critical and satisfies (\ref{assumptionweak}), and $\mathbf{P}_x[A (Y)>r]$ for any $x,r>0$.
If $A$ satisfies the monotonicity property (\ref{CBmonotonicity}),
then for any $\mathcal{F}_b$-measurable bounded random variable $F$,
as $r\rightarrow\infty$,
$$
\mathbf{E}_x[F| A(Y) >r ]\rightarrow \frac{1}{x}\mathbf{E}_x[Y_b\,F].
$$
\end{theo}

\begin{rem} \label{CBI}
For any $b\geq 0$ and $\mathcal{F}_b$-measurable bounded random variable $F$,
define a new probability $\mathbf{P}^*_x$
by
$$
\mathbf{E}^*_x[F]=\frac{1}{x}\mathbf{E}_x[Y_tF].
$$
It is well-known that $\mathbf{P}^*_x$ is the distribution of a CBI process $Y^*$ with the branching mechanism
$\Phi$ and the immigration mechanism $\Phi'$, where $\Phi'$ is the derivative of $\Phi$. See e.g., Section 2.3 and 3.1 in \cite{L12}
for some details on CBI processes.
Then Theorem \ref{LCCB} implies that for any $b>0$,
$(Y_t,t\in[0,b])$ under the conditioning of $\{A(Y)>r\}$ converges weakly to $(Y^*_t,t\in [0,b])$ as $r\rightarrow\infty$. In this case, we say that under the conditioning of $\{A(Y)>r\}$ the critical CB process $Y$ \emph{converges locally} to the CBI process $Y^*$.
\end{rem}

Now Theorem \ref{LCCB} and Lemma \ref{positive} immediately imply that the critical CB process $Y$ converges locally to $Y^*$ , under any of the three conditionings introduced in Section \ref{ss:L}.

\section{Continuum condensation trees}\label{s:CCT}

In this section first we define continuum condensation trees,
more precisely we define the left and right height processes of continuum condensation trees,
then we state two conjectures and one open problem on continuum condensation trees.

We now define continuum condensation trees.
Recall (\ref{Levy}) and (\ref{subordinator}).
Let $H$ be the height process associated with the L\'{e}vy process $X$ with the Laplace exponent $-\Phi$ and let $(H',X')$ be a copy of $(H,X)$.
Let $I=(I_t,t\geq 0)$ and $I'=(I'_t,t\geq 0)$ be the infimum processes of $X$ and $X'$ respectively.
Let $(U,V)$ be a bivariate subordinator with the Laplace exponent $\Phi(p,q)$.
Let $U^{-1}=(U^{-1}_t,t\geq 0)$ and $V^{-1}=(V^{-1}_t,t\geq 0)$ be the right-continuous inverses of $U$ and $V$ respectively.
Also introduce a random variable $\xi_\alpha$, which is exponential with parameter $\alpha$.
Assume that $(X,H)$, $(X',H')$, $(U,V)$,  and $\xi_\alpha$ are independent.
We define $\overleftarrow{H}'$ by
$$
\overleftarrow{H}'_t=H_t+U^{-1}_{-I_t}\quad \text{ if }U^{-1}_{-I_t}<\xi_\alpha, \quad \text{and} \quad
\overleftarrow{H}'_t=H_t+\xi_\alpha\quad\text{ if }U^{-1}_{-I_t}\geq\xi_\alpha, \quad t\geq 0.
$$
Then similarly define $\overrightarrow{H}'$ by
$$
\overrightarrow{H}'_t=H'_t+V^{-1}_{-I'_t}\quad \text{ if }V^{-1}_{-I'_t}<\xi_\alpha, \quad \text{and} \quad
\overrightarrow{H}'_t=H'_t+\xi_\alpha\quad\text{ if }V^{-1}_{-I'_t}\geq\xi_\alpha, \quad  t\geq 0.
$$
The processes $\overleftarrow{H}'$ and $\overrightarrow{H}'$ are called respectively left and right height processes of the continuum condensation tree with the branching mechanism $\Phi$.
From here on it is then natural to define the continuum condensation tree with the branching mechanism $\Phi$ as a random metric space $\mathcal{T}_*(\Phi)$ from the height processes $\overleftarrow{H}'$ and $\overrightarrow{H}'$.
For details refer to Page 103 in \cite{D08}.

Note that in the critical case, $\alpha=0$ implies that $\xi_\alpha=+\infty$ a.s., so the definitions of the height processes of continuum condensation trees and continuum immortal trees coincide.
Clearly the effect of $\xi_\alpha$ can also be achieved by killing $(U,V)$ at an independent exponential time with parameter $\alpha$, that is, we define
$$
\Phi'(p,q)=\frac{\Phi(p)-\Phi(q)}{p-q} \quad\text{for} \quad p\neq q, \quad\text{and}\quad \Phi'(p,p)=\Phi'(p),
$$
then let $(U',V')$ be a bivariate subordinator with the Laplace exponent $\Phi'(p,q)$. From here on we just proceed as in the definition of the left and right height processes of continuum immortal trees. However we feel that graphically speaking the above definition with $\xi_\alpha$ is somewhat clearer.

Then we turn to two conjectures and one open problem on continuum condensation trees. The first conjecture is related to the conditioning of large maximal degree, see Section \ref{ss:L} for an introduction of this conditioning. The maximal degree of L\'{e}vy trees corresponds to both the maximal out-degree of GW trees and the maximal jump of CB processes.
Recall that for GW trees, it is known that under the conditioning of large maximal out-degree, the local limit of a subcritical GW tree is a condensation tree. See Section 1 in \cite{AD14b} for a definition of condensation trees and see \cite{H14} for the proof of the local convergence to condensation trees.
It is also known that under the conditioning of large maximal jump, the local limit of a subcritical CB process with the branching mechanism $\Phi$ is a CBI process with the branching mechanism $\Phi$ and the immigration mechanism $\Phi'$, the derivative of $\Phi$. Note that $\Phi'(0)=\alpha>0$, so the CBI process is killed at an independent exponential time with parameter $\alpha$, where killing means sending to $\infty$.
For the proof see Theorem 4.4 in \cite{HL14}.
Recall our definition of continuum condensation trees, which is inspired by the definitions of condensation trees and the CBI processes that we just mentioned.
Then naturally we expect these continuum condensation trees to be the correct local limits of subcritical L\'{e}vy trees under the conditioning of large maximal degree. However, the desired proof seems to be much more involved than the two proofs we mentioned above, and currently we do not have it yet.
Now we state this expected local convergence explicitly as the following conjecture. Recall from Section \ref{ss:bp} that $\mathbf{N}$ is the excursion
measure of the strong Markov process $X-I$ at zero.
Also recall that we write $X$ for the canonical process under $\mathbf{N}$.

\begin{conj} \label{conjecturemd}
Assume that $\alpha>0$ and $\pi$ has unbounded support.
Then as $r\rightarrow\infty$,
$$
(H_{t\wedge \zeta},H_{(\zeta-t)_+};t\geq 0)
\text{ under }\mathbf{N}[\cdot|\sup \Delta X >r]\longrightarrow
(\overleftarrow{H}'_t,\overrightarrow{H}'_t;t\geq 0)
$$
weakly in $\mathbb{C}([0,\infty),\mathbb{R}^2)$.
\end{conj}

The second conjecture is related to the conditioning of large total mass, see Section \ref{ss:L} for an introduction of this conditioning. The total mass of L\'{e}vy trees corresponds to the total progeny of GW trees.
Recall that for GW trees, it is known that under the conditioning of large total progeny, the local limit of a subcritical GW tree
is an immortal tree or a condensation tree, depending on the offspring distribution. See \cite{J12,AD14b} for details.
In particular, if the subcritical offspring distribution $p=(p_0,p_1,p_2,\ldots)$ satisfies that $\sum_{k\in \mZ_+} a^k p_k=\infty$ for any $a>1$,
then the local limit of the subcritical GW tree $\tau(p)$ conditioned on large total progeny is the condensation tree with the offspring distribution $p$.
Inspired by this result, we make a conjecture on the local convergence of L\'{e}vy trees to continuum condensation trees, under the conditioning of large total mass.

\begin{conj}\label{conjecturetm}
Assume that $\alpha>0$ and for any $a>0$,
$$\int_1^\infty e^{a\theta}\pi (d\theta)=\infty.$$
Then as $r\rightarrow\infty$,
$$
(H_{t\wedge \zeta},H_{(\zeta-t)_+};t\geq 0)
\text{ under }\mathbf{N}[\cdot|\sigma > r]\longrightarrow
(\overleftarrow{H}'_t,\overrightarrow{H}'_t;t\geq 0)
$$
weakly in $\mathbb{C}([0,\infty),\mathbb{R}^2)$.
\end{conj}

We conclude this paper with an open problem on the conditioning of large width, see the end of Section \ref{ss:GWtail} for the definition of the width of GW trees and see also Section \ref{ss:L} for an introduction of this conditioning in the setting of L\'{e}vy trees.
Recall from Corollary \ref{GWwidth} and Corollary \ref{GWfour} that under the conditioning of large width, the local limit of a critical GW tree is an immortal tree.
See Corollary \ref{L3} for the corresponding result on L\'{e}vy trees. For subcritical GW trees and subcritical L\'{e}vy trees, the situation of the local convergence under this conditioning is completely unknown at the moment. In plain words, we want to know:

(Open Problem)
What is the local limit of a subcritical GW tree under the conditioning of large width?
What is the local limit of a subcritical L\'{e}vy tree under the conditioning of large width?

\bigskip
{\bf Acknowledgement}: Sincere thanks to an anonymous referee for many useful comments and suggestions, which improved considerably the presentation of this paper.


\end{document}